\documentclass[12pt]{article}

\newtheorem{theorem}{Theorem}

\newtheorem{corollary}[theorem]{Corollary}

\newtheorem{example}[theorem]{Example}

\newtheorem{lemma}[theorem]{Lemma}

\newtheorem{problem}[theorem]{Problem}

\newtheorem{remark}[theorem]{Remark}

\newenvironment{proof}[1][Proof]{\textbf{#1.} }{\ \rule{0.5em}{0.5em}}
\newdimen\dummy
\dummy=\oddsidemargin
\begin{document}
\date{}
\title{A monomorphism theorem for the inverse limit of nested retracts}
\author{Paul Fabel \\
Department of Mathematics \& Statistics\\
Mississippi State University}
\maketitle

\begin{abstract}
Suppose a given space $X$ can be realized as the inverse limit of nested
retracts $X_{n}$ where $X_{n}$ admits a universal cover or $X_{n}$ has a
certain more general property.

In various ways this paper characterizes injectivity of the canonical
homomorphism (into the inverse limit of the fundamental groups of the
factors). For example injectivity is equivalent to the existence of a kind
of generalized universal cover over $X$, a fibration whose fibres have
trivial path components and whose total space is simply connected.
\end{abstract}

\section{Introduction}

A paper of Biss \cite{Biss} shows \bigskip quite generally that the based
fundamental group of a topological space $X$ admits a canonical topology and
moreover the topological group $\pi _{1}(X,p)$ is invariant under the
homotopy type of $X$. If $X$ admits a universal cover in the usual sense
then $\pi _{1}(X,p)$ has the discrete topology. However, both the algebraic
and topological structure of $\pi _{1}(X,p)$ can be challenging to
understand if $X$ fails to be locally contractible. This paper aims to
explore $\pi _{1}(X,p)$ in the context of inverse limit spaces $%
X=\lim_{\leftarrow }X_{n},$ where $X$ can be approximated by factors $X_{n}$
where $\pi _{1}(X_{n},p)$ is either discrete or totally disconnected.

For example if $X=HE$ is the Hawaiian earring, \ (the one point union of a
null sequence of simple closed curves $X_{n}$ joined at a common point $p$)
the canonical homomorphism into the inverse limit of free groups $\phi :\pi
_{1}(HE,p)\rightarrow \lim_{\leftarrow }\pi _{1}(X_{n},p)$ is one to one
(Theorem 4.1 Morgan/Morrison \cite{mor}). The injectivity of $\phi $ ensures
that $\pi _{1}(HE,p)$ is a $T_{1}$ space, and since $\pi _{1}(HE,p)$ is a
topological group $\pi _{1}(HE,p)$ is completely regular. Further analysis
reveals that $\phi $ is not surjective, and it is shown in \cite{fab2} and 
\cite{fab5} that $\phi $ fails to be a homeomorphism onto its image, and
that $\pi _{1}(HE,p)$ fails to be a Baire space despite being the
uncountable regular topological fundamental group of a Peano continuum.

In general the homomorphism $\phi $ can fail to be one to one and $\pi
_{1}(X,p)$ can fail to be a $T_{1}$ space. For example the harmonic
archipelago $HA$ discussed in \cite{Biss} \cite{Bog} \cite{fab1} has an
uncountable topological fundamental group endowed with the course topology,
despite the fact that $HA$ is the inverse limit of nested absolute retracts.

The contrasting properties of the two examples $HE$ and $HA$ suggest the
possibility of a useful theory underlying the investigation of topological
fundamental groups.

For some positive results, it is shown in \cite{fab6} that the extra
topological structure of $\pi _{1}(X,p)$ can sometimes distinguish the
homotopy type of spaces with isomorphic homotopy groups, offsetting the
general failure of the Whitehead theorem. It is shown in \cite{fab4} that
the condition $\pi _{1}(X,p)$ is $T_{1}$ is equivalent to each retract of $X$
determining a closed topological fundamental subgroup. This in turn yields
`no retraction' theorems such as: The Hawaiian earring cannot be embedded as
a retract of any space whose topological fundamental group is completely
metrizable.

The above discussion highlights the role that injectivity of $\phi :\pi
_{1}(X)\rightarrow \lim_{\leftarrow }\pi _{1}(X_{n})$ can play in the
investigation of $\pi _{1}(X).$ This paper develops a list of
characterizations of injectivity of $\phi $ in terms of the topology of $\pi
_{1}(X,p)$ (Theorem \ref{main})$.$ As an application we obtain a
classification of the existence of a certain kind of generalized universal
covering space: the injectivity of $\phi $ characterizes the existence of a
fibration $q:X_{p}^{\symbol{126}}\rightarrow X$ such that $X_{p}^{\symbol{126%
}}$ is simply connected and each fibre is totally path disconnected
(Corollary \ref{main3}).

The results in this paper are related to the work of a number of other
authors including Cannon, Conner \cite{can}, de Smit \cite{desmit},Eda \cite
{eda}, Fischer and Zastrow \cite{Fischer}.

We close the paper with some related open questions in hopes of spurring
further developments in the algebraic topology of nonlocally contractible
spaces.

\section{Definitions and Preliminaries}

All of the following definitions are compatible with those found in Munkres 
\cite{Munkres}.

Suppose $X$ is a metrizable space and $p\in X.$ Let $S_{p}(X)=\{f:[0,1]%
\rightarrow X$ such that $f$ is continuous and $f(0)=p\}.$ Endow $S_{p}(X)$
with the topology of uniform convergence.

Let $C_{p}(X)=\{f\in S_{p}(X)|f(1)=p\}.$

The \textbf{topological fundamental group} $\pi _{1}(X,p)$ is the set of
path components of $C_{p}(X)$ endowed with the quotient topology under the
canonical surjection $q:C_{p}(X)\rightarrow \pi _{1}(X,p)$ satisfying $%
q(f)=q(g)$ if and only of $f$ and $g$ belong to the same path component of $%
C_{p}(X).$

Thus a set $U\subset \pi _{1}(X,p)$ is open in $\pi _{1}(X,p)$ if and only
if $q^{-1}(U)$ is open in $\pi _{1}(X,p).$

If $X\subset Y$ then $X$ is a \textbf{retract} of $Y$ if there exists a map $%
f:Y\rightarrow X$ such that $f_{X}=id_{X}.$ The space $X$ is $T_{1}$ if each
one point subset of $X$ is closed.

A space $X$ is \textbf{semilocally simply connected} at $p$ if $p\in X$ and
there exists an open set $U$ such that $p\in U$ and $j_{U}:U\rightarrow X$
induces the trivial homomorphism $j_{U}^{\ast }:\pi _{1}(U,p)\rightarrow \pi
_{1}(X,p).$

A topological space $X$ is \textbf{discrete} if every subset is both open
and closed. The space $X$ is \textbf{totally disconnected} if each the
components of $X$ are the one point subsets.

If $A_{1},A_{2},...$ are topological spaces and $f_{n}:A_{n+1}\rightarrow
A_{n}$ is a continuous surjection then, (endowing $A_{1}\times A_{2}..$ with
the product topology) the \textbf{inverse limit space }$\lim_{\leftarrow
}A_{n}=\{(a_{1},a_{2},...)\in (A_{1}\times A_{2}...)|f_{n}(a_{n+1})=a_{n}\}.$

A surjective map $q:E\rightarrow B$ \ is a \textbf{fibration} provided for
each space $Y$, for each map $f:Y\rightarrow E,$ and for each map $F:Y\times
\lbrack 0,1]\rightarrow B$ such that $F(y,0)=q(f)(0),$ there exists a map $%
F^{\symbol{126}}:Y\times \lbrack 0,1]\rightarrow E$ such that $q(F^{\symbol{%
126}})=F$ and $F^{\symbol{126}}(y,0)=q(y).$

\begin{remark}
The topological fundamental group $\pi _{1}(X,p)$ is a topological group
under concatenation of paths. (Proposition 3.1\cite{Biss}). A map $%
f:X\rightarrow Y$ determines a continuous homomorphism $f^{\ast }:\pi
_{1}(X,p)\rightarrow \pi _{1}(Y,f(p))$ via $f^{\ast }([\alpha ])=[f(\alpha
)] $ (Proposition 3.3 \cite{Biss}). If the spaces $X$ and $Y$ have the same
homotopy type then there is an isomorphism between $\pi _{1}(X,p)$ and $\pi
_{1}(Y,p)$ which is a homeomorphism.
\end{remark}

\section{Admissible inverse limit spaces}

The results in the paper apply to inverse limit spaces $X$ whose factors
have totally disconnected fundamental group. Suppose $\{p_{1}\}\subset
X_{1}\subset X_{2}...$ is a nested sequence of path connected separable
metric spaces and suppose $r_{n}:X_{n+1}\rightarrow X_{n}$ is a retraction
and suppose $\pi _{1}(X_{n},p)$ is totally disconnected. For the remainder
of this paper the resulting inverse limit space $X=\lim_{\leftarrow }X_{n}$
is said to be \textbf{admissible. }Given such an admissible space $X$ we
also define $p=(p_{1},p_{1},...)\in X.$ Let $R_{n}:X\rightarrow X_{n}$
denote the canonical map $R_{n}(x_{1},x_{2},...)=x_{n}.$ Let $R_{n}^{\ast
}:\pi _{1}(X,p)\rightarrow \pi _{1}(X_{n},p_{1})$ denote the induced
homomorphism. Let $\phi :\pi _{1}(X,p)\rightarrow \lim_{\leftarrow }\pi
_{1}(X_{n},p_{n})$ denote the induced homomorphism $\phi ([f])=(R_{1}^{\ast
}([f]),R_{2}^{\ast }([f]),...).$ Let $P$ denote the constant map of $%
C_{p}(X).$ Let $[P]$ denote the path component of $P$ in $C_{p}(X).$

\begin{remark}
If $X_{n}$ is both locally path connected and semilocally simply connected (
for example if $X_{n}$ is an $ANR$ or if $X_{n}$ is locally contractible)
then $\pi _{1}(X)$ is discrete and in particular $\pi _{1}(X_{n})$ is
totally disconnected. \cite{fab2}.
\end{remark}

\begin{lemma}
\label{ezlem}The map $j_{n}:X_{n}\rightarrow \lim_{\leftarrow }X_{n}$
defined via $j_{n}(x_{n})=(x_{1},..x_{n},x_{n},..)$ is an embedding of $%
X_{n} $ onto a retract of $X.$ Henceforth we can treat $X_{n}$ as a retract
of $X.$
\end{lemma}

\begin{proof}
Note the function $j_{n}$ is well defined since $r_{n+k}$ fixes $X_{n}$
pointwise. Moreover $j_{n}$ is continuous since $\Pi _{k}(j_{n})$ is
continuous with $\Pi _{k}:\lim_{\leftarrow }X_{n}\rightarrow X_{k}$ the
projection map. The map $j_{n}$ is one to one since if $x_{n}\neq y_{n}$
then $\Pi _{n}(x_{n})\neq \Pi _{n}(y_{n}).$ To check that $j_{n}$ is a
homeomorphism onto its image note if $\lim_{k\rightarrow \infty }$ $%
(x_{1}^{k},...x_{n}^{k},x_{n}^{k},..)=(x_{1},...x_{n},x_{n},...)$ then
convergence is coordinatewise and in particular $x_{n}^{k}\rightarrow x_{n}.$
Finally the map $R_{n}^{\symbol{94}}:\lim_{\leftarrow }X_{n}\rightarrow
im(j_{n})$ defined via $R_{n}^{\symbol{94}}=(\Pi _{1},\Pi _{2},..\Pi
_{n},\Pi _{n},..)$ is the desired retraction.
\end{proof}

\section{A monomorphism theorem}

For admissible spaces $X$ Theorem \ref{main} offers various
characterizations of injectivity of the canonical homomorphism $\phi :\pi
_{1}(X,p)\rightarrow \lim_{\leftarrow }\pi _{1}(X_{n},p).$ It is interesting
to note that whether or not $\phi $ is injective is completely determined by
topological properties of $\pi _{1}(X,p).$

\begin{theorem}
\label{main}Suppose $X$ is admissible as demonstrated by the retractions $%
r_{n}:X_{n+1}\rightarrow X_{n}$ and suppose $\phi :\pi _{1}(X,p)\rightarrow
\lim_{\leftarrow }\pi _{1}(X_{n},p)$ is the canonical homomorphism. The
following are equivalent.
\end{theorem}

\begin{enumerate}
\item  $\phi $ is one to one.

\item  $\pi _{1}(X,p)$ is $T_{1}$.

\item  $\pi _{1}(X,p)$ is normal.

\item  $\pi _{1}(X,p)$ is totally disconnected.

\item  $\pi _{1}(X,p)$ is totally path disconnected.
\end{enumerate}

\begin{proof}
First we establish the equivalence of $1,2,$ and $3.$

$1\Rightarrow 2.$ Suppose $\phi $ is one to one. Since $\pi _{1}(X_{n},p)$
is $T_{1},$ the product $\Pi _{n=1}^{\infty }\pi _{1}(X_{n},p)$ is $T_{1}.$
Since $\phi $ is continuous and one to one $\phi ^{-1}(\phi (x))=\{x\}$ is
closed. Thus $\pi _{1}(X,p)$ is $T_{1}.$

$2\Rightarrow 3.$ Suppose $\pi _{1}(X,p)$ is $T_{1}.$ Since $\pi _{1}(X,p)$
is a topological group and since $\pi _{1}(X,p)$ is $T_{1},$ $\pi _{1}(X,p)$
is regular (ex 6 p.145 \cite{Munkres}). Since $X$ is a separable metric
space $X,$ the Urysohn metrization theorem shows $X$ can be embedded in the
Hilbert cube (Thm 4.1 p. 217 \cite{Munkres}). Hence $C_{p}(X)$ is a
separable metric space and in particular $C_{p}(X)$ is Lindelof. Thus $\pi
_{1}(X,p)$ is Lindelof since the quotient map $q:C_{p}(X)\rightarrow \pi
_{1}(X,p)$ is surjective. Hence $\pi _{1}(X,p)$ is paracompact since $\pi
_{1}(X,p)$ is both regular and Lindelof (Ex 2 p. 259 \cite{Munkres}). Thus,
since $\pi _{1}(X,p)$ is paracompact, $\pi _{1}(X,p)$ is normal (Thm. 4.1 p.
255 \cite{Munkres}).

$3\Rightarrow 2$ by definition$.$

To show $2\Rightarrow 1$ suppose $\pi _{1}(X,p)$ is $T_{1}$ and $[f]\in \ker
(\phi ).$ Let $f=(f_{1},f_{2},..).$ The retracts $r_{n+k}:X_{n+k+1}%
\rightarrow X_{n+k}$ enable us to define $f^{n}\in C_{p}(X)$ via $%
f^{n}=(f_{1},..f_{n},f_{n},..).$ Lemma \ref{ezlem} shows $%
j_{n}(R_{n}):X\rightarrow im(j_{n})$ is a retraction such that $%
j_{n}(R_{n})(f)=f^{n}.$ Suppose $[f]\in \ker (\phi ).$ Since $[f]\in \ker
(\phi ),$ $f^{n}$ is inessential in $j_{n}(X_{n})$ for each $n$ and in
particular $f^{n}$ is inessential in $X.$ Note $f^{n}\rightarrow f$ in $%
C_{p}(X).$ Since $\pi _{1}(X)$ is a $T_{1}$ space the trivial element is
closed in $\pi _{1}(X,p)$. Since $q:C_{p}(X)\rightarrow \pi _{1}(X,p)$ is a
quotient map the path component $[P]$ is closed in $C_{p}(X).$ Since $%
f^{n}\rightarrow f$ in $C_{p}(X)$ and since $f^{n}$ is path homotopic to $P$
it follows that $f$ is path homotopic to $P.$ Hence $\phi $ is one to one.

Now we prove the equivalence of $2,4,$ and $5.$

$4\Rightarrow 5$ trivially.

$5\Rightarrow 1.$ Suppose $5$ holds. Let $\overline{\{e\}}$ denote the
closure of the trivial element $\{e\}$ in $\pi _{1}(X,p).$ Notice every
function $\alpha :[0,1]\rightarrow \overline{\{e\}}$ is continuous. Thus if $%
\pi _{1}(X,p)$ fails to be $T_{1}$ then $\overline{\{P\}}\neq \{P\}$ and
hence $5$ would fail to hold since selecting $f\in \overline{\{P\}}%
\backslash \{P\}$ be may define a (continuous) function $\alpha
:[0,1]\rightarrow \{e,f\}$ such that $\alpha (t)=e$ if $0\leq t<1$ and $%
\alpha (1)=f.$ Thus $5\Rightarrow 1.$

$1\Rightarrow 4.$ If $\pi _{1}(X,p)$ is $T_{1}$ then $\phi $ is continuous
and one to one and maps $\pi _{1}(X,p)$ into the totally disconnected space $%
\Pi \pi _{1}(X_{n},p).$ Hence $\pi _{1}(X,p)$ is totally disconnected. $%
4\Rightarrow 5$ trivially.

This completes the proof.
\end{proof}

\section{Application to generalized universal covers}

Corollary \ref{main3} offers another characterization of injectivity of $%
\phi $ in terms of the existence of a kind of generalized universal cover
over $X.$

The paper \cite{Biss} develops the following generalization of the familiar
covering maps.

A fibration $q:E\rightarrow X$ is a \textbf{rigid covering fibration}
provided the following 3 properties hold.

\begin{enumerate}
\item  For each $x\in X$ the fibre $q^{-1}(x)$ is totally path disconnected.

\item  If $n\neq 1$ then the induced homomorphism $q_{n}^{\ast }:\pi
_{n}(E)\rightarrow \pi _{1}(X)$ is an isomorphism.

\item  If $n\geq 2$ then the induced homomorphism $q_{1}^{\ast }:\pi
_{1}(E)\rightarrow \pi _{1}(X)$ is one to one.
\end{enumerate}

\begin{remark}
It shown in Spanier \cite{Spanier} that a fibration has unique path lifting
if and only if each fibre is totally path disconnected (Lemma 4 p.68). For
such fibrations (Theorem 4 p.72, Corollary 11 p.377) $q_{1}^{\ast }$ is one
to one and $q_{n}^{\ast }$ is an isomorphism for $n\geq 2.$ Consequently
conditions $2$ and $3$ are consequences of condition 1.
\end{remark}

The familiar universal cover can be seen as a special case of a \textbf{%
generalized universal cover, }a rigid fibration whose total space is simply
connected. Locally path connected spaces which are semilocally simply
connected spaces admit a universal cover. Which path connected spaces admit
a generalized universal cover? Theorem \ref{main2} provides a complete
answer. We should point out that Theorem \ref{main2} also follows by
combining Theorems 4.3 and 4.6 and Corollary 4.7 of \cite{Biss}. For the
sake of clarity we include proof of Theorem \ref{main2}.

\begin{theorem}
\label{main2} Suppose $X$ is any path connected space. Then $\pi _{1}(X,p)$
is totally path disconnected if and only if there exists a rigid covering
fibration $q:E\rightarrow X$ such that $E$ is simply connected.
\end{theorem}

\begin{proof}
Suppose $\pi _{1}(X,p)$ is totally path disconnected. Let $E$ denote the
following quotient space of $S_{p}(X).$ Declare $f\symbol{126}g$ if $f$ and $%
g$ are path homotopic in $X.$ Define $Q:E\rightarrow X$ such that $%
Q(f)=f(1). $ To prove $Q:E\rightarrow X$ is a fibration suppose $%
f:Y\rightarrow E$ is any map and $F:Y\times \lbrack 0,1]\rightarrow X$ is
any homotopy such that $F(y,0)=Q(f(y)).$ For each $y\in Y$ select $\beta
_{y}\in f(y).$ For each $y\in Y$ and $t\in \lbrack 0,1]$ define $\alpha
_{y,t}:[0,1]\rightarrow X$ such that $\alpha _{y,t}(s)=F(y,st).$ Define $F^{%
\symbol{126}}:Y\times \lbrack 0,1]\rightarrow E$ such that $F^{\symbol{126}%
}(y,t)=[\beta _{y}\ast \alpha _{y,t}]$ where $\ast $ denotes the familiar
concatenation of paths. To check that $F^{\symbol{126}}$ is well defined
note $\beta _{y}(1)=F(y,0)=\alpha _{y,0}(s).$ Since $\alpha _{y,0}$ is
constant $\beta _{y}\ast \alpha _{y,0}$ is path homotopic to $\beta _{y}.$
Thus $F^{\symbol{126}}(y,0)=f(y).$ To check that $F^{\symbol{126}}$ is
continuous let $A=\{(f,g)\in S_{p}(X)\times S_{p}(X)|f(1)=g(0)\}.$ Let $%
B=\{[f],[g]\in E\times E|(f,g)\in A\}.$ Note $(Q\times Q):A\rightarrow B$ is
a quotient map. Since we are using the uniform topology on $S_{p}(X)$ the
standard path concatenation function $\ast :A\rightarrow E$ is continuous.
To check that $\ast $ induces a map from $B\rightarrow E$ it suffices to
note that $f\ast g$ is path homotopic to $f^{^{\prime }}\ast g^{^{\prime }}$
if $[f]=[f^{^{\prime }}]$ and $[g]=[g^{^{\prime }}].$ Thus $Q:E\rightarrow X$
is a fibration. By definition $Q^{-1}(p)$ is canonically homeomorphic to the
topological fundamental group $\pi _{1}(X,p).$ If $x\in X$ and $\alpha $ is
a path from $p$ to $x$ then the fibre over $x$ is homeomorphic to $\pi
_{1}(X,p)$ via $\pi _{1}(X,p)\ast \lbrack \alpha ].$ Thus $Q$ is a rigid
covering fibration with unique path lifting. To prove that $E$ is simply
connected note if $\alpha :[0,1]\rightarrow X$ is an inessential loop then
the (unique) lift of $\alpha $ based at $[P]$ is the function $\alpha ^{%
\symbol{126}}:[0,1]\rightarrow E$ satisfying $\alpha ^{\symbol{126}%
}(t)=[\beta ^{t}]$ and $\beta ^{t}:[0,1]\rightarrow X$ satisfies $\beta
^{t}(s)=\alpha (st)$. In particular $\alpha ^{\symbol{126}}(1)=[\alpha ]\neq
\lbrack P].$ Thus if $\gamma :[0,1]\rightarrow E$ is a loop based at $P$
then $Q(\gamma )$ is inessential. Hence $Q(\gamma )$ bounds a disk and the
(unique) lift of the disk bounding $Q(\gamma )$ bounds $\gamma .$ Thus $E$
is simply connected.

Conversely, suppose $q:E\rightarrow X$ is a rigid covering fibration such
that $E$ is simply connected. Consider the following homotopy $%
F:C_{p}(X)\times \lbrack 0,1]\rightarrow X.$ Let $F(\alpha ,t)=\alpha (t).$
Select $e\in E$ such that $q(e)=p.$ Let $F^{\symbol{126}}:C_{p}(X)\times
\lbrack 0,1]\rightarrow E$ denote the unique lift satisfying $F^{\symbol{126}%
}(\alpha ,0)=e.$

Since $E$ is simply connected $F^{\symbol{126}}(\alpha ,1)=F^{\symbol{126}%
}(\beta ,1)$ if and only if $\alpha $ and $\beta $ are path homotopic in $X.$
Since $\pi _{1}(X,p)$ inherits the quotient topology, $F^{\symbol{126}}$
induces a continuous injection $g:\pi _{1}(X,p)\rightarrow q^{-1}(p)$
determined the rule $g([\alpha ])=F^{\symbol{126}}(\alpha ,1).$ Since $im(g)$
is totally path disconnected, and since $g$ is one to one, $\pi _{1}(X,p)$
must totally path disconnected.
\end{proof}

Combining Theorems \ref{main} and \ref{main2}, in the context of admissible
spaces, the existence of a simply connected rigid covering fibration is
characterized by injectivity of the canonical homomorphism.

\begin{corollary}
\label{main3}Suppose the inverse limit space $X$ is admissible as
demonstrated by the retracts $r_{n}:X_{n+1}\rightarrow X_{n}$ and suppose $%
\phi :\pi _{1}(X,p)\rightarrow \lim_{\leftarrow }\pi _{1}(X_{n},p)$ is the
canonical homomorphism. Then the following are equivalent:
\end{corollary}

\begin{enumerate}
\item  $\phi $ is one to one.

\item  There exists a rigid covering fibration $q:E\rightarrow X$ such that $%
E$ is simply connected.
\end{enumerate}

\begin{example}
The Hawaiian earring $HE$ is locally path connected but not semilocally
simply connected. Consequently $\pi _{1}(HE)$ fails to be discrete (Thm 2 
\cite{fab3}) and $HE$ does not admit a universal cover in the usual sense.
However the canonical homomorphism $\phi :\pi _{1}(HE)\rightarrow
\lim_{\leftarrow }\pi _{1}(F_{n})$ is one to one (Thm. 4.1 \cite{mor}).
Consequently $HE$ admits a generalized universal cover.
\end{example}

\section{Questions}

Here are some natural questions relevant to the results of this paper.

\begin{problem}
Which path connected separable metric spaces $X$ are admissible?
\end{problem}

\begin{problem}
Suppose $X$ is any path connected space such that $\pi _{1}(X,p)$ is $T_{1}.$
Must $\pi _{1}(X,p)$ be totally disconnected?
\end{problem}

\begin{problem}
Theorem 4.8 \cite{Biss} asserts that if $q:E\rightarrow X$ is a rigid
covering fibration such that $E$ is locally path connected and simply
connected then $\pi _{1}(X,p)$ is isomorphic to the group of homeomorphisms $%
f:E\rightarrow E$ such that $q(f)=q.$ To what extent can the assumption that 
$E$ is locally path connected be dropped? Under what conditions is the
isomorphism also a homeomorphism with respect to the compact open topology?
\end{problem}

\end{document}